# Dyck Numbers, IV. Nested patterns in OEIS A036991

Gennady Eremin

***Abstract.*** Dyck paths (also balanced brackets and Dyck words) are among the most heavily studied Catalan families. This paper is a continuation of [2, 3, 4]. In the paper we are dealing with the numbering of Dyck paths, with the resulting numbers, the terms of the OEIS sequence A036991, which encode Dyck paths and which we have called the Dyck numbers. In this sequence, it turned out to be promising to investigate nested term patterns. The sequence is constructed entirely from nested patterns, with each pattern having an infinite number of copies. We have already considered repeated triplets of adjacent odd numbers; copies of such triplets in the process of expansion form quite complex structures, the analysis of which allows us to understand the structure of A036991. Additionally, we note that the composition of the terms of each pattern is repeated many times by a simple shift (offset) of Dyck numbers by some fixed value. As a result, pattern copies give us important parts of most levels and even give us new levels.

## 1    Introduction

Dyck paths (balanced brackets or Dyck words) are among the most heavily studied Catalan families [1]. This article continues [2, 3, 4]. In the paper we are dealing with the numbering of Dyck paths, with the resulting numbers, the terms of the OEIS A036991, which encode Dyck paths and which we have called the *Dyck numbers*. In the previous papers, we looked at some features of the OEIS sequences A036991 and A350346. This paper examines nested data patterns.

We have already considered *triplets* of adjacent odd numbers of the form $t$-4, $t$-2 and $t$, where $t$, for example, can be the $n$th Mersenne number (see the OEIS A000225); but there are also an infinite number of triples of such adjacent triplets, aka *nine-terms*. There are many copies of larger patterns that allow you to understand the structure of A036991. Suffice it to say that the list of terms of each level of the sequence is repeated many times by simply shifting the Dyck numbers by some fixed value. As a result, such copies give us essential parts of subsequent levels and even give us new levels.

For example, we can consider as a pattern all the terms with a certain length of binary code. A simple copy of such a pattern is obtained by adding a binary unit at the beginning of each term, which simply changes the length of the binary code and increases the selected terms by a fixed value (the distance between adjacent terms does not change). This procedure can be repeated many times, and as a result we will get an arbitrary number of copies of our pattern.

In the sequence A036991 (as well as A350346), the terms are ordered in ascending order; let's show the initial A036991 terms $a_n$, $n \geq 1$ (the initial term 0 corresponds to the empty Dyck path, Mersenne numbers are shown in red):

0, 1, 3, 5, 7, 11, 13, 15, 19, 21, 23, 27, 29, 31, 39, 43, 45, 47, 51, 53, 55, 59, 61, 63, 71, 75, 77, 79, 83, 85, 87, 91, 93, 95, 103, 107, 109, 111, 115, 117, 119, 123, 125, 127, 143, 151, 155, 157...

The group of terms with the same binary length $n$ forms the *n-level*, $Ł_n$, and this is clearly visible in A350346. The $n$-level contains terms from a half-open interval that is bounded by two Mersenne numbers $(M_{n-1}; M_n]$, where $M_n = 2^n - 1$, $n \geq 1$. Therefore, each subsequent level contains natural numbers from the double length segment. In this regard, the number of terms in the next level is often doubled. The maximum of *n*-level is always the *n*th Mersenne number (in binary we have *n* ones, no zeros), and it can be called a *level top*. The *n*-level minimum is determined by the *Dyck successor function*, *DSucc-function*, for $M_{n-1}$ [2].

In the sequence A036991 the size of the levels (the number of Dyck numbers) is determined by the [A001405](#) terms $a_n$, $n \geq 1$: 1, 1, 2, 3, 6, 10, 20, 35, 70, 126, 252, 462… Let's write it this way $\#Ł_n = A001405(n-1)$, $n > 0$. When moving from an even level to an odd level, the number of terms doubles. For example, at 4-level we have the first triplet (11,13,15), and level 5 already has two triplets (19,21,23) and (27,29,31).

In almost all levels of the sequence A036991, the initial terms are far from obvious; however, we know well the final term of each level, which is the corresponding Mersenne number. It is also known that the final term, starting at level 4, is a senior element of the triplet. Next we will show that at the following levels the Mersenne number is a senior term of very large patterns that are repeatedly nested within each other. Moreover, there are often repeated nested patterns, and it is more convenient to view such patterns from the end of the level. Therefore, it is logical to analyze each level in A036991 from the end, using a known *Predecessor function*.

## 2 Predecessor function.

In mathematics, the 1-ary Predecessor function *Pred* maps a natural number to the previous natural number, that is, $Pred(n) = n - 1$, $n > 0$. Additionally, $Pred(0) = 0$. In our case we are working with Dyck number, so for a given A036991 term *t* the *Dyck Predecessor function*, *DPred-function*, should return the preceding A036991 term, denote $DPred(t)$. For example, $DPred(39) = 31$. It is logical to assume that the predecessor function is undefined in the initial term 0.

Next, we will need the predecessor function mainly to work not with individual terms, but with a group (tuple) of ordered consecutive (adjacent) A036991 terms. (Although in the general case, we can consider a single term as a tuple of length 1.)

**Definition 1.** Let's call a *pattern* a set (tuple) of ascending adjacent A036991 terms of the same binary length *n*: $x = (x_1, …, x_k)$, provided that

$$M_{n-1} < x_1 = DPred(x_2), x_2 = DPred(x_3), …, x_{k-1} = DPred(x_k), x_k \leq M_n.$$

The length of *x* is $len(x) = x_k - DPred(x_1)$, the power is $\#x = k$. Obviously, $k \leq \#Ł_n$.

Often we will work with patterns that include all terms of a certain level. For example, in level 4 there are three terms 11, 13, and 15 that form the first triplet; denote such a pattern $\pi_4(15) = \pi_4(M_4) = (11, 13, 15) = Ł_4$. In parentheses, we will specify the senior term of the pattern, *pattern top* (this is often the Mersenne number). Obviously, the length of the pattern is $len(\pi_4(M_4)) = M_4 - M_3 = 8$ and $\#\pi_4(M_4) = 3$. Further the pattern $\pi_4$, like other patterns, is repeated many times, let's say, *copied*.

**Definition 2.** Let the pattern $x = (x_1, …, x_k)$ be given. The tuple $y = (y_1, …, y_k)$ is called a *copy* of *x* if

(i) $x_k < y_1$ and

(ii) $DPred(y_1) - DPred(x_1) = y_1 - x_1 = y_2 - x_2 = \ldots = y_k - x_k = \Delta$.

This equality may not hold for subsequent terms, i.e. perhaps $DSucc(y_k) - DSucc(x_k) \neq \Delta$. Let's call $\Delta$ the *offset* (shift) between the pattern and its copy and denote by *offset*$(x, y)$.

Usually the maxima of patterns and copies are available to us (often these are Mersenne numbers), so the offset between the pattern $x$ and its copy $y$ is easier to define as follows: *offset*$(x, y) = \max y - \max x$. In the same way we will fix the offset between two copies of the same pattern.

We have already noted that if you add one at the beginning of the binary code of each term of an arbitrary pattern, you get a copy at the next level with the same number of terms and the same distance between adjacent terms. This procedure can be repeated endlessly, going over levels. The following statement is obvious.

**Proposition 3.** Each pattern has an infinite number of copies.

In copies, as in patterns, we also specify senior terms. For example, here are the final triplets for some levels: $\pi_4(M_5) = (M_5 - 4, M_5 - 2, M_5)$, $\pi_4(M_6) = (59, 61, \mathbf{63})$, $\pi_4(511) = (M_9 - 4, M_9 - 2, M_9)$. The senior term in parentheses (the *triplet top*) varies in copies and is always indicated. But in the case of level 4, where the pattern $\pi_4$ is generated, we can omit the brackets, for example, we can write so $\pi_4 = (11, 13, \mathbf{15})$. Obviously, for all copies $len(\pi_4) = 8$, $\#\pi_4 = 3$. We will reduce in the general case: $\pi_n = \pi_n(M_n) = (DSucc(M_{n-1}), \ldots, M_n)$, $\#\pi_n = \#L_n$. In the following, we will often mark the senior term of the tuple (pattern) in bold.

**Definition 4.** Let the pattern $x = (x_1, \ldots, x_k)$ and its copy $y = (y_1, \ldots, y_k)$ be given. Let us say that *x precedes y*, that is, $x = DPred(y)$ if $x_k = DPred(y_1)$.

Often triplets are directly adjacent to each other. Here's a small example.

**Example 5.** In level 5 we have triplets $\pi_4(23) = (19, 21, \mathbf{23})$ and $(27, 29, \mathbf{31}) = \pi_4(31)$. Both triplets are adjacent because $23 = DPred(27)$. In this case $\pi_4(23) = DPred(\pi_4(31))$, and the offset of such copies is *offset*$(\pi_4(23), \pi_4(31)) = 31 - 23 = 8 = len(\pi_4)$. □

# 3 Structural patterns

## 3.1. **Odd levels.**

In the odd level 5, we have two copies $\pi_4(M_5 - 8)$ and $\pi_4(M_5)$. The doubling of the pattern $\pi_4$ is due to the transition from odd suffixes of length 3, 3-suffixes, to even suffixes of length 4, 4-suffixes (see [4], proposition 2). We have three valid odd 3-suffixes 011, 101, and 111 with positive dynamics (the number of ones exceeds the number of zeros); recall, the binary code of an odd number always ends with unit, and suffix 001 is not valid. Obviously, when expanding, we can add either a 0 or a 1 at the beginning of each 3-suffix. As a result, we get six 4-suffixes and respectively two copies of the pattern $\pi_4$ at 5-level (as before, in each copy of $\pi_4$ we have three terms). Let's show the 5-level as a connection (union) of two copies of pattern $\pi_4$:

$$L_5 = \pi_4(M_5 - 8) \oplus \pi_4(M_5)$$
$$= DPred(\pi_4(M_5)) \oplus \pi_4(M_5) = \pi_4(M_5)^2 = (19, 21, 23, 27, 29, \mathbf{31}).$$

We can consider the resulting expression as the next pattern $\pi_5(M_5)$ or $\pi_5$ for short, which also has copies (with different senior terms); but such copies are in little demand. Maybe this is because there are no original (unusual) components in the odd-level patterns. Obviously, $\#\pi_5 = 2 \times \#\pi_4 = 6$.

For compact notation, we *doubled* the pattern $\pi_4$ using the degree $2$; this doubling resembles the concatenation operation, only the shift of Dyck numbers by the value of $len(\pi_4)$ must be taken into account. In this case, the copy $\pi_4(M_5)$ is docked to the preceding copy $\pi_4(M_5 - len(\pi_4))$. What is interesting, at the next level we have three adjacent triplets, and there we get a *triple copy* (using the degree $3$):

$$\pi_4(M_6 - 2 \times len(\pi_4)) \oplus \pi_4(M_6 - len(\pi_4)) \oplus \pi_4(M_6) = \pi_4(M_6)^3 \quad \text{or}$$
$$\pi_4(47) \oplus \pi_4(55) \oplus \pi_4(63) = \pi_4(63)^3 = (43, 45, 47, 51, 53, 55, 59, 61, \mathbf{63}).$$

Obviously, in the general case

$$\pi_n(t)^2 = DPred(\pi_n(t)) \oplus \pi_n(t) = \pi_n(t - len(\pi_n)) \oplus \pi_n(t), \quad t > M_n, \text{ and}$$
$$\pi_n(t)^3 = DPred(DPred(\pi_n(t))) \oplus DPred(\pi_n(t)) \oplus \pi_n(t), \quad t - len(\pi_n) > M_n.$$

Each pattern, as well as its copy, is placed on some segment of natural numbers. Let us give a more general definition of the length of an arbitrary pattern.

**Definition 6.** The *length* of the pattern $\pi_n$ is $len(\pi_n) = M_n - M_{n-1} = 2^{n-1}$.

Obviously, $DPred(\pi_n(t)) = \pi_n(t - 2^{n-1})$ for the senior term $t > M_{n-1}$. Accordingly, $len(\pi_4) = M_4 - M_3 = 8$, $len(\pi_5) = M_5 - M_4 = 16$, and so on. It is important not to confuse pattern length $len(\pi_n)$ (the interval on the numeric axis) and its power $\#\pi_n$ (the number of A036991 terms on this interval). Let us formulate the obvious statement.

**Proposition 7.** For an odd $n$-level, $n \geq 5$,

$$\pi_n = \pi_n(M_n) = DPred(\pi_{n-1}(M_n)) \oplus \pi_{n-1}(M_n)$$
$$= \pi_{n-1}(M_n - 2^{n-2}) \oplus \pi_{n-1}(M_n)$$
$$= \pi_{n-1}(M_n)^2. \tag{1}$$

As you can see, with odd levels everything is quite simple: *the odd level pattern is a doubled copy of the preceding even level pattern*.

When we move from an even level to the next odd level, doubling the number of terms is consistent with doubling the copies of the corresponding pattern. At level 5 we have two copies of the pattern $\pi_4$ with senior terms $31 = M_5$ and $23 = (M_5 + M_4)/2 = M_5 - 8$. Term 23 is the center of 5-level. Starting at level 5, such central terms (in both odd and even levels) form the sequence $\{a_n\}_{n \geq 5} = \{23, 47, 95, 191, 383, 767, 1535,...\}$. Let's denote $H_n = (M_{n-1} + M_n)/2 = M_{n-1} + 2^{n-2} = M_n - 2^{n-2}$, $n \geq 5$; obviously, $H_{n+1} = 2H_n + 1$ (the recurrence relation is identical to the Mersenne numbers, and this is natural). We can say that each central term $H_n$ divides the $n$-level (the interval $M_n - M_{n-1} = 2^{n-1}$) into two equal segments $(M_{n-1}; H_n]$ and $(H_n; M_n]$. The segment on the high side of $M_n$ always contains a copy of the pattern of $(n-1)$-level, $\pi_{n-1}(M_n)$; in this case the offset is $2^{n-1}$.

In the OEIS, there are several sequences with similar central terms that differ only in the initial numbers and indexing. For example, you can choose A052940 or A290114, where all the terms are Dyck numbers. There are few formulas in these sequences, and we can add the following formulas (taking into account the accepted notation in the OEIS and the indexing of terms):

$$A052940(n) = A000225(n+1) + 2\wedge n = A000225(n+2) - 2\wedge n$$
$$= (A000225(n+1) + A000225(n+2))/2, \text{ where } n \geq 1;$$
$$A290114(n) = A000225(n) + 2\wedge(n-1) = A000225(n+1) - 2\wedge(n-1)$$
$$= (A000225(n) + A000225(n+1))/2, \text{ where } n > 1.$$

Next, let's look at the 6th level.

### 3.2. **Even levels.**

3.2.1. At level 5 we have 6 terms, two triplets $\pi_4(23)$ and $\pi_4(31)$. In the binary code of each such term there are 5 bits; recall that in all Dyck numbers except zero, the highest bit and the zero (right-sided) bit are fixed and equal to 1. To make the transition from 5-level to 6-level, we need to extend the binary 4-suffixes in the terms. We can add unit at the beginning of each 4-suffix (the suffix dynamics do not deteriorate), as a result we get 6 terms of 6-level. Obviously, such a procedure simply increases the numbers by $2^5 = 32$, so we get the last two triplets of the 6th level: $\pi_4(55)$ and $\pi_4(63)$. In general, for the even $n$-level, $n \geq 6$, on the side of maximum $M_n$ at the interval $(M_n - M_{n-1})/2 = 2^{n-2}$, there are two copies of the pattern of $(n-2)$-level:

$$DPred(\pi_{n-2}(M_n)) \oplus \pi_{n-2}(M_n) = \pi_{n-2}(M_n)^{\underline{2}}.$$

We get additional central terms (in both even and odd levels) $a_n$, $n \geq 6$: 55, 111, 223, 447, 895, 1791... with equalities $a_n = (H_n + M_n)/2 = H_n + 2^{n-3} = M_n - 2^{n-3}$ and $a_{n+1} = 2a_n + 1$ (again the recurrence relation is identical to the Mersenne numbers). In the OEIS, there is the corresponding sequence [A086224](#) (we ignore the initial terms). There are few formulas in this sequence, and we can add the following formula (taking into account the accepted notation and indexing of terms):

$$A086224(n) = A052940(n+1) + 2\wedge n = A000225(n+3) - 2\wedge n$$
$$= (A052940(n+1) + A000225(n+3))/2, \text{ where } n \geq 0.$$

3.2.2. At level 6 we have another third triplet, the copy $\pi_4(47) = \pi_4(H_6) = (43, 45, \mathbf{47})$. We get the third triplet if in the copy $\pi_4(31)$, making a transition from level 5 to level 6, in the terms at the beginning of each suffix of length 4 (this is the maximum suffix in the terms of level 5) we add not a binary unit, but zero. But we need to show that the dynamics of the new suffixes will not become negative. Indeed, when we went from the pattern $\pi_4(15)$ to copy $\pi_4(31)$, to each suffix of length 3 (such 3-suffixes have a balance $\geq 1$), we added a binary unit at the beginning. Therefore, in terms of $\pi_4(31)$ the 4-suffixes have a balance $\geq 2$. The further transition from $\pi_4(31)$ to $\pi_4(47)$ is made by adding a leading zero to the 4-suffixes, so the new 5-suffixes have the same balance as the original 3-suffixes. In fact, in pattern $\pi_4(15)$, we added binary fragment 01 to each 3-suffix, as a result, the dynamics of the new 5-suffixes does not change.

  Thus in level 6 we have three triplets: the triplet in the middle $\pi_4(47)$, the subsequent triplet $\pi_4(55)$ and the top triplet $\pi_4(63)$. The last two triplets are a copy of the 5th level pattern $\pi_5(63) = \pi_4(55) \oplus \pi_4(63) = \pi_4(63)^{\underline{2}}$. We can also write $\pi_4(H_6) \oplus \pi_5(M_6) = \pi_4(M_6)^{\underline{3}}$.

Recall that each triplet has three terms, therefore we have already received 9 terms of the 6th level (by the way, there are 10 terms in level 6). Let us formulate another statement.

**Proposition 8.** Each even $n$-level, $n \geq 6$, ends with tuple

$$\pi_{n-2}(H_n) \oplus \pi_{n-1}(M_n) = \pi_{n-2}(H_n) \oplus \pi_{n-2}(M_n)^2$$
$$= \pi_{n-2}(M_n)^3. \qquad (2)$$

Recall that in this case $\pi_{n-2}(H_n) = DPred(DPred(\pi_{n-2}(M_n)))$. Obviously, in (2) the number of terms is $3 \times \#\pi_{n-2}$.

Thus, each even $n$-level, $n \geq 6$, is divided into 4 equal segments; the upper three segments (starting from the upper boundary $M_n$) contain three copies of pattern $\pi_{n-2}$. The length of each such segment is $(M_n - M_{n-1})/4 = 2^{n-3}$. The terms of the three copies of the pattern $\pi_{n-2}$ are computed by a simple fixed shift of the pattern terms. Serious problems arise when calculating the terms of the last segment (these are the initial terms of the $n$-level, near the lower bound $M_{n-1}$), where we need to get some semblance of another copy of pattern $\pi_{n-2}$. We will deal with these problems in the rest of this paper.

3.2.3. As you can see, at even levels, everything is more complicated. As in the case of equality (1), we also start from an even level, obtaining most of the terms of the next even level, but the initial part (less than a quarter) of the new level terms is not yet achievable. This can be clearly seen in the analysis of binary suffixes in terms. In general, to obtain terms of the even $n$-level, we need to analyze in ($n$-2)-level the odd binary ($n$-3)-suffixes in the terms and increase such suffixes by two digits if possible. For this we use the obvious four *binary fragments* 00, 01, 10, and 11. Let's look at these fragments in reverse order.

*Fragment 11 (F11).* If in each term of ($n$-2)-level we add the fragment F11 at the beginning of the maximal suffix (length $n-3$), we get the top copy $\pi_{n-2}(M_n)$ at level $n$. In the terms of this copy, the balance of suffixes increased by 2 (compared to the balance of maximal suffixes in the original pattern $\pi_{n-2}$). Adding F11 gives us an increase (shift) in the terms of the resulting copy by $M_n - M_{n-2} = 4 \times 2^{n-2} - 2^{n-2} = 3 \times 2^{n-2}$.

*Fragment 10 (F10).* Adding F10 gives us the next copy $DPred(\pi_{n-2}(M_n)) = \pi_{n-2}(M_n - 2^{n-3})$ which precedes the copy discussed above. In the terms of the second copy, the balance of maximal suffixes does not change (coincides with the balance of suffixes in the original pattern $\pi_{n-2}$). The terms of the second copy are increased by $3 \times 2^{n-2} - 2^{n-3} = 5 \times 2^{n-3}$.

*Fragment 01 (F01).* The third fragment F01 gives us the copy $DPred(DPred(\pi_{n-2}(M_n))) = \pi_{n-2}(M_n - 2^{n-2})$; the senior term of this third copy is located in the center of level $n$. In the terms of the third copy, the balance of maximal suffixes does not change too. The terms of the third copy are increased by $5 \times 2^{n-3} - 2^{n-3} = 4 \times 2^{n-3} = 2^{n-1}$.

*Fragment 00 (F00).* The last fragment F00 significantly decreases the dynamics of maximal suffixes in the terms of the original pattern $\pi_{n-2}$, so some terms are rejected and do not appear at $n$-level. For example, in 4-level, only suffix 111 retains the valid dynamics after adding two leading zeros; and only this suffix generates the term in 6-level $39 = 100111_2$. The remaining two suffixes 011 and 101 turn into Dyck words 0011 and 0101 when the first leading zero is added (the dynamic is 0 in a Dyck word), so another leading zero is not allowed. The natural question arises: how many even $k$-level terms are rejected by the fragment F00 in the general case? The answer is obvious, as many as there are Dyck words of length $k$ (equivalently, we can get a Dyck word of length $k$ if we simply zero out the highest binary bit in the term of an even $k$-level). Thus, the fragment F00 gives us a truncated copy of pattern $\pi_{n-2}$, and for this fragment we can reject in advance terms of the original pattern $\pi_{n-2}$, those terms in which the dynamics is 1 (the difference between the number of ones and zeros in the binary code).

3.2.4. We can count the number of $n$-level terms that have passed the "sieve" of the fragment F00. As we know, Dyck words are counted by Catalan numbers: Cat($k$) is the number of Dyck words of length $2k$ (see the OEIS A000108 in [5]). Therefore, the number of terms of an even $n$-level that are rejected by the fragment F00 is Cat($n/2 - 1$). The terms of the shortened copy are placed in the lower segment of the $n$-level, i.e. they begin an even $n$-level. Let us denote the shortened copy by $\mu_n$. Obviously,

$$\#\mu_n = \#\pi_{n-2} - \text{Cat}(n/2 - 1). \tag{3}$$

**Definition 9.** The shortened copy of an even level has a special status; let's call the tuple $\mu_n$ (the beginning of $n$-level) an *n-level core* or *n-core* for short.

Formula (3) gives us a sequence of numbers that enumerate A036991 terms in cores in even $n$-levels, $n = 6, 8, 10, \ldots$ Here is the beginning of this sequence:

1, 5, 21, 84, 330, 1287, 5005, 19448, 75582, 293930, 1144066, 4457400...

This is a known sequence of binomial coefficients of the following form (see the OEIS A002054):

$$C(2k+1, k-1) = k \times \text{Cat}(k+1)/2, \text{ where } k \geq 1. \tag{4}$$

But this sequence does not yet have a formula similar to (3). Taking into account the OFFSET of terms in sequences and the notation adopted in the OEIS, formula (3) takes the following form:

$$\text{A002054}(k) = \text{A001405}(2k+1) - \text{A000108}(k+1), \ k \geq 1. \tag{5}$$

We can now formulate the obvious statement.

**Proposition 10.** For an even $n$-level, $n \geq 6$,

$$\#\pi_n = 4 \times \#\pi_{n-2} - \text{Cat}(n/2 - 1) = 2 \times \#\pi_{n-1} - \text{Cat}(n/2 - 1).$$

In the OEIS A001405, there is a similar formula by Christopher Hanusa (2003).

It is easy to see that at an even $n$-level the third (full) copy of pattern $\pi_{n-2}$ and the core $\mu_n$ are placed in the segment $(M_{n-1}; H_n]$. The full copy $\pi_{n-2}(H_n)$ is on the right side of this segment. The core $\mu_n$ ends with the senior term $H_n - len(\pi_{n-2}) = H_n - 2^{n-3} = M_{n-1} + 2^{n-3}$. For example, at 6-level in $\mu_6$ there is a single term $H_6 - len(\pi_4) = 47 - 8 = 39$, i.e. $\mu_6 = \mu_6(39) = (39)$. Then at level 8 the senior term of $\mu_8$ is $H_8 - len(\pi_6) = 191 - 32 = 159$, and so on. As a result, we get additional central terms (in both even and odd levels) $a_n$, $n \geq 6$: 39, 79, 159, 319, 639, 1279, 2559... Again the recurrence relation is identical to the Mersenne numbers, $a_{n+1} = 2a_n + 1$.

In the OEIS, there is the corresponding sequence A052549 (and we again ignore the initial terms). There are few formulas in this sequence, and we can add the following formula (taking into account the accepted notation and indexing of terms):

$$A052549(n) = A000225(n+1) + 2^{\wedge}(n-1), \quad n \geq 0.$$

## 4 Generation of levels in A036991

The analysis of patterns and their copies allows us to understand the structure of the sequence A036991. Each odd level of the sequence is easily generated from the preceding even level. The odd $(n+1)$-level contains two copies of all the terms of the even $n$-level, and in each copy the terms increase by a fixed value: in the upper copy (the senior term is $M_{n+1}$) we increase the even-level terms by $M_{n+1} - M_n = 2^n$, in the next copy the even-level terms increase by only half, i.e. the offset is $2^{n-1}$. Calculation of terms in even levels is much more complicated and more interesting (in this regard, there is an understandable desire to simplify the sequence A036991 by removing all odd levels, because any odd level can be easily reconstructed if necessary).

Table 1 below contains the data of the initial levels (the first three levels are obvious and therefore not shown; the reader can easily extend the table if necessary). The last line is added for orientation; the numbers of this line are not used in the calculations (so it is shown pale). The cells of the table contain the senior terms of patterns and their copies. The red color shows the senior terms in the cores. The rows begin with the terms that are directly used in the calculations. For example, line $M_n$ begins with the senior term 15 of the first triplet, pattern $\pi_4(M_4) = (11,13,\mathbf{15})$ in column 4 (even 4-level). The line $H_n$ begins with the senior term 23 of the second copy $\pi_4(H_5) = \pi_4(23) = (19,21,\mathbf{23})$ in column 5 (odd level).

| Level $n$ | 4 | 5 | 6 | 7 | 8 | 9 | 10 | 11 | 12 | 13 | 14 | 15 |
|---|---|---|---|---|---|---|---|---|---|---|---|---|
| $M_n$ | 15 | 31 | 63 | 127 | 255 | 511 | 1023 | 2047 | 4095 | 8191 | 16383 | 32767 |
| $M_n - 2^{n-3}$ | | | 55 | | 223 | | 895 | | 3583 | | 14335 | |
| $H_n$ | | 23 | 47 | 95 | 191 | 383 | 767 | 1535 | 3071 | 6143 | 12287 | 24575 |
| $H_n - 2^{n-3}$ | | | 39 | | 159 | | 639 | | 2559 | | 10239 | |
| $M_{n-1}$ | 7 | 15 | 31 | 63 | 127 | 255 | 511 | 1023 | 2047 | 4095 | 8191 | 16383 |

Table 1. Senior terms of patterns and their copies.

In the odd columns, we have two senior terms, and in each even column (starting from the 6th level) we have four senior terms according to the breakdown into four segments in levels. Even-level terms are easier to calculate using previous even-level information. Let's look at an example.

**Example 11.** Consider that we know all ten terms of the even 6th level:

$\pi_6 = (39) \oplus \pi_4(47) \oplus \pi_4(55) \oplus \pi_4(63)$
  $= (39,\ 43,45,47,\ 51,53,55,\ 59,61,\mathbf{63}) = \pi_6(63).$

The length of such a nested pattern is $len(\pi_6) = M_6 - M_5 = 63 - 31 = 32$. Let's calculate the terms of the next even 8-level. First, we write three copies of the pattern $\pi_6$ with senior terms $255 - 2 \times 32 = 191$, $255 - 32 = 223$ and $255$; the terms of the pattern $\pi_6$ are increased in the copies by the following values, respectively: $191 - 63 = 128$, $223 - 63 = 160$ and $255 - 63 = 192$. The result is the next tripled copy:

$$\begin{aligned}\pi_6(M_8)^3 &= \pi_6(H_8) \oplus \pi_6(M_8 - len(\pi_6)) \oplus \pi_6(M_8) \\ &= \pi_6(191) \oplus \pi_6(223) \oplus \pi_6(255) \\ &= (167,\ 171,173,175,\ 179,181,183,\ 187,189,\mathbf{191}) \oplus \\ &\quad (199,\ 203,205,207,\ 211,213,215,\ 219,221,\mathbf{223}) \oplus \\ &\quad (231,\ 235,237,239,\ 243,245,247,\ 251,253,\mathbf{255}).\end{aligned}$$

The senior term of the 8-core is $255 - 3 \times 32 = 127 + 32 = 159$. In order to get the core $\mu_8 = \mu_8(159)$, we will analyze the binary codes. Let's show the binary codes of the ten terms of level 6 (see A350346):

<u>100111</u>, <u>101011</u>, <u>101101</u>, 101111, <u>110011</u>, <u>110101</u>, 110111, 111011, 111101, 111111.

We have underlined five terms. Such terms are rejected in $\mu_8$ because the fragment F00 produces a negative dynamics of the maximal suffixes in these terms (in binary code, the number of 0's exceeds the number of 1's). Let's add that such terms generate Dyck words, which are counted by Catalan numbers. We increase the remaining five terms 47, 55, 59, 61, and 63 by $159 - 63 = 96$. As a result, we obtain the initial terms of level 8: $\mu_8(159) = (143,151,\ 155,157,159) = (143,151) \oplus \pi_4(159)$. Thus, all 35 terms of 8th level are obtained. We can write $\#\mu_8 = 5$ and $\#\pi_8 = 5 + 3 \times \#\pi_6 = 35$. □

The last example shows how, starting from an arbitrary even $n$-level, it is quite simple to calculate the terms of the subsequent odd $(n+1)$-level and even $(n+2)$-level. However, some problems arise in obtaining the initial terms of the $(n+2)$-level, the core of the level $n+2$. In Example 11, we had to analyze the binary codes of ten terms of the initial level 6, and as a result we selected five terms with which we obtained the first terms of the 8th level, the core $\mu_8(159)$. And here a natural question arises, is it possible to simplify the procedure for selecting terms for $\mu_8(159)$? How to avoid analyzing binary expansions?

## 5 Cores in even levels

Dyck numbers are of two kinds [3], they are (i) *triplet terms* and (ii) *ternary tree roots*. The root of a ternary tree is generally distant from adjacent terms at a distance of 4 or more (the exception is root 1). In other words, a given A036991 term $r > 1$ is a root if both $r-2$ and $r+2$ are not Dyck numbers. Every $n$th Mersenne number, $n > 4$, generates tree root using the *Dyck Successor function* (see [2], Theorem 5):

$$DSucc(M_n) = M_n + 2^m, \text{ where } m = \lceil n/2 \rceil.$$

For example, the term 39 mentioned above (the first term in level 6) is the tree root, which is generated by the 5th Mersenne number, i.e. $DSucc(M_5) = 31 + 2^3 = 39$. On subsequent levels, each Mersenne number generates a group of roots. For example, at the even 10-level we get the following initial terms, or roots of ternary trees: $DSucc(M_9) =$

543, $DSucc(543) = 559$, and $DSucc(559) = 567$. Next we have the triplet $\pi_4(575) =$ (571, 573, **575**) at the level 10 (and this is not all the terms of the 10-core).

Now let's go back to Example 11 and repeat the initial terms of 8-level:

$$\mu_8(159) = (143, 151, 155, 157, 159) = (143, 151) \oplus \pi_4(159).$$

Obviously, $len(\mu_8(159)) = 159 - DPred(143) = 159 - 127 = 32$. As we see, in $\mu_8(159)$ the first two terms are the roots of the ternary tree, which are generated by the 7th Mersenne number: $DSucc(M_7) = 143$ and $DSucc(DSucc(M_7)) = 151$.

In the 8-core, the two roots are followed by the triplet $\pi_4(159)$. Recall that in the sequence A036991, each term $t$ from $n$-level (it can be both a tree root and a triplet term) generates one triplet of the form $\pi_4(4t + 3) = (4t - 1, 4t + 1, 4t + 3)$ at $(n+2)$-level. And vice versa, each triplet of $(n+2)$-level has a single parent in the preceding $n$-level. In our case, the senior term of the triplet is $4t + 3 = 159$, so the triplet $\pi_4(159)$ is generated by the term $t = (159 - 3)/4 = 39$, and it is the only term of the 6-core.

Note that term 39 generates a triplet for $\mu_8(159)$, but the term 39 itself is rejected in the 8-core, since the fragment F00 gives negative suffix dynamics. And this is true for all the roots in each core. Each root from the core of an even $n$-level is rejected by the fragment F00 in the $(n+2)$-core.

If we look at the senior terms in Table 1, we find that each senior term of the $n$-level generates a triplet with a senior term in a similar segment of the $(n+2)$-level. For example, the $n$th Mersenne number generates a triplet

$$(4M_n - 1, 4M_n + 1, 4M_n + 3) = (M_{n+2} - 4, M_{n+2} - 2, M_{n+2}).$$

And the senior term $H_n$ generates a similar triplet

$$(4H_n - 1, 4H_n + 1, 4H_n + 3) = (H_{n+2} - 4, H_{n+2} - 2, H_{n+2}) \text{ and so on.}$$

This is directly related to the recurrence relation:

$$a_{n+2} = 2a_{n+1} + 1 = 2(2a_n + 1) + 1 = 4a_n + 3.$$

Thus, we can formulate the following statement.

**Proposition 12.** An arbitrary A036991 term $t$ from some fixed segment of the even $n$-level generates the triplet $(4t - 1, 4t + 1, 4t + 3)$ in the same segment of the $(n+2)$-level.

The considered procedure for generating triplets allows us to obtain a significant fraction of terms in even-level cores, since each term of an arbitrary core is guaranteed to generate three terms in the next core. Here's another example.

**Example 13**. From the data of the core $\mu_8(159)$, we can easily obtain all the triplets of the next core $\mu_{10}(4 \times 159 + 3) = \mu_{10}(639)$. Five terms of the 8-core generate 5 triplets in 10-core, as a result, we immediately get 15 terms in the next 10-core:

(571, 573, 575), (603, 605, 607), (619, 621, 623), (627, 629, 631), (635, 637, 639).

According to formula (3) the number of terms in the 10-core is

$$\#\mu_{10} = \#\pi_8 - \text{Cat}(10/2 - 1) = 35 - 14 = 21.$$

Obviously, the remaining 6 terms are tree roots. The first three roots are generated by the Mersenne number $M_9 = 511$: $DSucc(511) = 543$, $DSucc(543) = 559$ and $DSucc(559) = 567$. The first root is shifted by $2^{\lceil 9/2 \rceil} = 2^5 = 32$ from the 9th Mersenne number, then the shift is halved in steps, it is 16 and 8. The remaining three roots 591, 599 and 615 are not achievable for us yet. □

## 6 Core subsequence of A036991, core patterns

In the sequence A036991, the distance between adjacent terms is a power of 2, namely 2, 4, 8, 16 and so on (this is directly indicated in the corresponding comments to A036991 in [5]). The length of the levels in A036991 is also a power of 2, the difference between two adjacent Mersenne numbers. The length of the $n$-level, aka the length of the pattern $π_n$, is equal to $len(π_n) = M_n - M_{n-1} = 2^{n-1}$.

Starting at level 8, we begin to divide each level into 4 equal segments, and this allowed us to better navigate the structure of the levels (the length of the segment is also a power of 2, $2^{n-1}/4 = 2^{n-3}$). For example, at the even level $n+2$, the pattern $π_n$ is repeated in the top three segments, starting from the senior term $M_{n+2}$ (with a simple shift of the pattern terms). But it is more complicated with the cores in even levels (in Example 13 we failed to calculate all the roots in the 10-core). Let's select cores of even levels (starting with 6-level) and compose the corresponding subsequence of A036991:

39, 143, 151, 155, 157, 159, 543, 559, 567, 571, 573, 575, 591, 599, 603, 605, 607, 615, 619, 621, 623, 627, 629, 631, 635, 637, 639, 2111, 2143, 2159, 2167, 2171, 2173,…

The senior terms of the cores are shown in red. There are 500 initial terms of the core subsequence in the Appendix below.

Next, we will begin split the cores, each core also divided into four equal intervals, four *subsegments*. Such fragmentation is possible if the segment length is sufficient (at $n$-level, the length of the subsegment is $2^{n-3}/4 = 2^{n-5}$). At level 10, we can already work with subsegments; let's continue with Example 13.

**Example 14** (*continuation of Example 13*). At level 10, the core $μ_{10}(639)$ is divided into 4 subsegments of length $2^{10-5} = (639 - M_9)/4 = 32$. Let's list the senior terms of the subsegments (shown in red): $M_9 + 32 = 543, 575, 607$, and $607+32 = 639$. Next, we show all core terms line by line (one line is one subsegment):

$μ_{10}$ = (543) ⊕
(559,567) ⊕ (571,573,575) ⊕
(591,599) ⊕ (603,605,607) ⊕
(615) ⊕ (619,621,623) ⊕ (627,629,631) ⊕ (635,637,639).

In such a record, all six roots of the core are clearly visible. In the first line we have a single term and it is the first root of the core $μ_{10}$. The second and third lines are identical, here we have two copies of the 8-level core $μ_8 = μ_8(159)$: $μ_8(575)$ and $μ_8(607)$. In this case, the core $μ_8$ manifests itself as a *core pattern* that includes two roots and one triplet. In these two lines we get 4 more roots. And finally, the last fourth line is a copy of the pattern $π_6$, which includes one root (the 6th root of the core) and three triplets, i.e. (615) ⊕ $π_4(639)^3 = π_6(639)$. Ultimately, we can briefly write down the 10-core as follows:

$μ_{10}$ = (543) ⊕ $μ_8(607)^2$ ⊕ $π_6(639)$. □

It may be noted that we have practically already begun to work with the core patterns $\mu_8 = \mu_8(159)$ and $\mu_{10} = \mu_{10}(639)$, and this is only the beginning. So far we have five patterns: two patterns for levels (level patterns) and three patterns for cores (core patterns). Let's repeat these patterns.

$\pi_4 = \pi_4(15) = (11, 13, 15)$, $len(\pi_4) = 8$, $\#\pi_4 = 3$;
$\pi_6 = \pi_6(63) = (39) \oplus \pi_4(47) \oplus \pi_4(55) \oplus \pi_4(63) = (39) \oplus \pi_4(63)^3$, $len(\pi_6) = 32$, $\#\pi_4 = 10$;
$\mu_6 = \mu_6(39) = (39)$, $len(\mu_6) = 8$, $\#\mu_6 = 1$;
$\mu_8 = \mu_8(159) = (143, 151) \oplus \pi_4(159)$, $len(\mu_8) = 32$, $\#\mu_8 = 5$;
$\mu_{10} = \mu_{10}(639) = (543) \oplus \mu_8(575) \oplus \mu_8(607) \oplus \pi_6(639)$, $len(\mu_{80}) = 128$, $\#\mu_{10} = 21$.

In the patterns, the senior terms are shown in red. There is only one term in the core pattern $\mu_6$, so such a pattern is practically useless and therefore rarely used. The pattern $\mu_{10}$ is of sufficient length, and we divided the 10-core into four subsegments. We have shown in blue the senior terms of the subsegments; the last term 639 is both the senior term of the top subsegment and the senior term of the 10-core.

Note that in pattern $\mu_{10}$, the second and third subsegments are copies of pattern $\mu_8$. Let's check this kind of copying at the next 12th core, consider another example.

**Example 15.** In the Appendix, the reader will find 84 terms of the 12-core, from the first term 2111 to the last term 2559 (the senior of the core is shown in red). The blue color shows the terms 2175, 2303 and 2431 (the senior terms of the subsegments). Packing the terms into patterns gives us the following compact expressions (the reader can easily check):

1st subsegment $(2111, 2143) \oplus \mu_8(2175) = \mu_{12/1}(2175)$;
2nd subsegment $(2207) \oplus \mu_8(2239) \oplus \mu_8(2271) \oplus \pi_6(2303) = \mu_{10}(2303)$;
3rd subsegment $(2335) \oplus \mu_8(2367) \oplus \mu_8(2399) \oplus \pi_6(2431) = \mu_{10}(2431)$;
4th subsegment $\mu_8(2463) \oplus \pi_6(2495) \oplus \pi_6(2527) \oplus \pi_6(2559) = \mu_{12/4}(2559)$.

We introduced two new core patterns $\mu_{12/1}$ and $\mu_{12/4}$, respectively, for the initial subsegment and the last (top) subsegment of the 12-core. In the following cores these patterns are actively used. As a result, we get

$\mu_{12} = \mu_{12}(2559) = \mu_{12/1}(2175) \oplus \mu_{10}(2303) \oplus \mu_{10}(2431) \oplus \mu_{12/4}(2559)$, $\#\mu_{12} = 84$.

As we see, the two internal subsegments are again copies of the previous core $\mu_{10}$. □

The Appendix has another 330 terms of the next 14th core. Let us show only the final results of the analysis (the reader can do his own calculations):

$\mu_{14} = \mu_{14}(10239) = \mu_{14/1}(8703) \oplus \mu_{12}(9215) \oplus \mu_{12}(9727) \oplus \mu_{14/4}(10239)$, $\#\mu_{14} = 330$.

where $\mu_{14/1}(8703) = (8319) \oplus \mu_{12/1}(8447) \oplus \mu_{12/1}(8575) \oplus \mu_{10}(8703)$
and $\mu_{14/4}(10239) = \mu_{10}(9855) \oplus \mu_{12/4}(9983) \oplus \mu_{12/4}(10111) \oplus \mu_{12/4}(10239)$.

Again, the preceding core pattern $\mu_{12}$ is completely copied in the second and third subsegments of the pattern $\mu_{14}$. While we cannot prove the necessary theorem, let us formulate a conjecture.

**Conjecture 16.** The core pattern $\mu_n$ ($n = 8, 10, 12,\ldots$) is copied in the second and third subsegments of the pattern $\mu_{n+2}$.

Apparently, to prove Hypothesis 16, we will have to work again with binary codes of the A036991 terms.

Let us calculate the offset of the core terms of an even $n$-level, $n \geq 8$, to obtain the second and third subsegments in the pattern $\mu_{n+2}$. The senior term of the $n$-core is $M_{n-1} + 2^{n-3}$, respectively, the senior term of the $(n+2)$-core is $M_{n+1} + 2^{n-1}$. Next, we need to descend in the pattern $\mu_{n+2} = \mu_{n+2}(M_{n+1} + 2^{n-1})$ by one subsegment and then one more subsegment (the length of the subsegment in $\mu_{n+2}$ is $2^{n-1}/4 = 2^{n-3}$). As a result, we get the following two offsets to get the required copies:

$$(M_{n+1} + 2^{n-1}) - (M_{n-1} + 2^{n-3}) - 2^{n-3} = 7 \times 2^{n-2} \quad \text{and} \quad 7 \times 2^{n-2} - 2^{n-3} = 13 \times 2^{n-3}. \tag{6}$$

Here's a little example.

**Example 17.** Let's consider the case of $n = 20$:

the senior term of the 20-core is $M_{19} + 2^{17} = 524287 + 131072 = 655359$,
the senior term of the 22-core is $M_{21} + 2^{19} = 2097151 + 524288 = 2621439$,
the length of the subsegment in the 22-core is $2^{17} = 131072$,
the senior term of the subsegment $\mu_{22/3}$ is $2621439 - 131072 = 2490367$,
the senior term of the subsegment $\mu_{22/2}$ is $2490367 - 131072 = 2359295$,
the offset for the terms of the subsegment $\mu_{22/3}$ is $2490367 - 655359 = 1835008 = 7 \times 2^{18}$,
the offset for the terms of the subsegment $\mu_{22/2}$ is $2359295 - 655359 = 1703936 = 13 \times 2^{17}$.

The obtained result confirms the formula (6). □

## 7 Conclusion

The proof of Conjecture 17 is probably a matter of time. But we do not yet know the initial and the last (top) subsegments in the level cores; what is interesting is the structure of these subsegments. Let's take a look at the top subsegments, starting with the 10-core:

$\mu_{10/4} = \pi_6(639)$, $\#\mu_{10/4} = 10$;
$\mu_{12/4} = \mu_8(2463) \oplus \mu_{10/4}(2495) \oplus \mu_{10/4}(2527) \oplus \mu_{10/4}(2559)$, $\#\mu_{12/4} = 35$;
$\mu_{14/4} = \mu_{10}(9855) \oplus \mu_{12/4}(9983) \oplus \mu_{12/4}(10111) \oplus \mu_{12/4}(10239)$, $\#\mu_{14/4} = 21 + 3 \times 35 = 126$;
$\mu_{16/4} = \mu_{12}(39423) \oplus \mu_{14/4}(39935) \oplus \mu_{14/4}(40447) \oplus \mu_{14/4}(40959)$, $\#\mu_{16/4} = 84 + 3 \times 126 = 462$.

As you can see, the terms of the top subsegments are counted by certain numbers of the OEIS sequence A001405. We can formulate another hypothesis.

**Conjecture 18.** The top subsegment of the $n$-core ($n = 12, 14, 16, ...$) is

$\mu_{n/4} = \mu_{n-4}(A - 3 \times 2^{n-7}) \oplus \mu_{(n-2)/4}(A)^3$, where $A = M_{n-1} + 2^{n-3}$ is the senior term of the $n$-core.

With the initial subsegment $\mu_{n/1}$ it turned out to be more complicated. We can say that $\mu_{n/1}$ is a new core inside the $n$-core, a *secondary core*. It can be assumed that within such a secondary core we will also find another core (this process is apparently infinite).

**Acknowledgement.** The author would like to thank Olena G. Kachko (Kharkiv National University of Radio Electronics) for helpful discussions.

# References

[1] R. Stanley. *Catalan numbers*. Cambridge University Press, Cambridge, 2015.

[2] G. Eremin. *Dyck numbers, I. Successor function*, 2022. arXiv preprint arXiv:2210.00744

[3] G. Eremin. *Dyck numbers, II. Triplets and rooted trees in OEIS A036991*, 2022. arXiv preprint arXiv:2211.01135

[4] G. Eremin. *Dyck numbers, III. Enumeration and bijection with symmetric Dyck paths*, 2023. arXiv preprint arXiv:2302.02765

[5] N. J. A. Sloane, The On-Line Encyclopedia of Integer Sequences, https://oeis.org.

Concerned with sequences A000108, A000225, A001405, A002054, A036991, A052549, A052940, A086224, A290114, and A350346.

*Email address*: ergenns@gmail.com
Written: June 17, 2023

# Appendix. Core subsequence of A036991 (the senior terms are shown in red).

39, 143, 151, 155, 157, 159, 543, 559, 567, 571, 573, 575, 591, 599, 603, 605, 607, 615, 619, 621, 623, 627, 629, 631, 635, 637, 639, 2111, 2143, 2159, 2167, 2171, 2173, 2175, 2207, 2223, 2231, 2235, 2237, 2239, 2255, 2263, 2267, 2269, 2271, 2279, 2283, 2285, 2287, 2291, 2293, 2295, 2299, 2301, 2303, 2335, 2351, 2359, 2363, 2365, 2367, 2383, 2391, 2395, 2397, 2399, 2407, 2411, 2413, 2415, 2419, 2421, 2423, 2427, 2429, 2431, 2447, 2455, 2459, 2461, 2463, 2471, 2475, 2477, 2479, 2483, 2485, 2487, 2491, 2493, 2495, 2503, 2507, 2509, 2511, 2515, 2517, 2519, 2523, 2525, 2527, 2535, 2539, 2541, 2543, 2547, 2549, 2551, 2555, 2557, 2559, 8319, 8383, 8415, 8431, 8439, 8443, 8445, 8447, 8511, 8543, 8559, 8567, 8571, 8573, 8575, 8607, 8623, 8631, 8635, 8637, 8639, 8655, 8663, 8667, 8669, 8671, 8679, 8683, 8685, 8687, 8691, 8693, 8695, 8699, 8701, 8703, 8767, 8799, 8815, 8823, 8827, 8829, 8831, 8863, 8879, 8887, 8891, 8893, 8895, 8911, 8919, 8923, 8925, 8927, 8935, 8939, 8941, 8943, 8947, 8949, 8951, 8955, 8957, 8959, 8991, 9007, 9015, 9019, 9021, 9023, 9039, 9047, 9051, 9053, 9055, 9063, 9067, 9069, 9071, 9075, 9077, 9079, 9083, 9085, 9087, 9103, 9111, 9115, 9117, 9119, 9127, 9131, 9133, 9135, 9139, 9141, 9143, 9147, 9149, 9151, 9159, 9163, 9165, 9167, 9171, 9173, 9175, 9179, 9181, 9183, 9191, 9195, 9197, 9199, 9203, 9205, 9207, 9211, 9213, 9215, 9279, 9311, 9327, 9335, 9339, 9341, 9343, 9375, 9391, 9399, 9403, 9405, 9407, 9423, 9431, 9435, 9437, 9439, 9447, 9451, 9453, 9455, 9459, 9461, 9463, 9467, 9469, 9471, 9503, 9519, 9527, 9531, 9533, 9535, 9551, 9559, 9563, 9565, 9567, 9575, 9579, 9581, 9583, 9587, 9589, 9591, 9595, 9597, 9599, 9615, 9623, 9627, 9629, 9631, 9639, 9643, 9645, 9647, 9651, 9653, 9655, 9659, 9661, 9663, 9671, 9675, 9677, 9679, 9683, 9685, 9687, 9691, 9693, 9695, 9703, 9707, 9709, 9711, 9715, 9717, 9719, 9723, 9725, 9727, 9759, 9775, 9783, 9787, 9789, 9791, 9807, 9815, 9819, 9821, 9823, 9831, 9835, 9837, 9839, 9843, 9845, 9847, 9851, 9853, 9855, 9871, 9879, 9883, 9885, 9887, 9895, 9899, 9901, 9903, 9907, 9909, 9911, 9915, 9917, 9919, 9927, 9931, 9933, 9935, 9939, 9941, 9943, 9947, 9949, 9951, 9959, 9963, 9965, 9967, 9971, 9973, 9975, 9979, 9981, 9983, 9999, 10007, 10011, 10013, 10015, 10023, 10027, 10029, 10031, 10035, 10037, 10039, 10043, 10045, 10047, 10055, 10059, 10061, 10063, 10067, 10069, 10071, 10075, 10077, 10079, 10087, 10091, 10093, 10095, 10099, 10101, 10103, 10107, 10109, 10111, 10127, 10135, 10139, 10141, 10143, 10151, 10155, 10157, 10159, 10163, 10165, 10167, 10171, 10173, 10175, 10183, 10187, 10189, 10191, 10195, 10197, 10199, 10203, 10205, 10207, 10215,

10219, 10221, 10223, 10227, 10229, 10231, 10235, 10237, 10239, 33023, 33151, 33215, 33247, 33263, 33271, 33275, 33277, 33279, 33407, 33471, 33503, 33519, 33527, 33531, 33533, 33535, 33599, 33631, 33647, 33655, 33659, 33661, 33663, 33695, 33711, 33719, 33723, 33725, 33727, 33743, 33751, 33755, 33757, 33759, 33767, 33771, 33773, 33775, 33779, 33781, 33783, 33787, 33789, 33791, 33919, 33983, 34015, 34031, 34039, 34043, 34045, 34047, 34111, 34143, 34159, 34167, 34171, 34173, 34175, 34207, 34223, 34231, 34235, 34237, 34239, 34255, 34263, 34267, 34269, 34271, 34279, 34283, 34285, 34287, 34291, 34293, 34295, 34299, 34301, 34303, 34367, 34399, 34415, 34423, 34427, 34429, 34431, 34463, 34479, 34487, 34491, 34493, 34495, 34511, 34519, 34523, 34525, 34527, 34535, 34539, 34541, 34543, 34547, 34549, 34551, 34555, 34557, 34559, 34591, 34607, 34615, 34619, 34621, 34623, 34639, 34647, 34651, 34653, 34655, 34663, 34667, 34669, 34671, 34675, 34677, 34679, 34683, 34685, 34687, 34703, 34711, 34715, 34717, 34719, 34727, 34731, 34733, 34735, 34739, 34741, 34743, 34747, 34749, 34751, 34759, 34763, 34765, 34767, 34771, 34773, 34775, 34779, 34781, 34783, 34791, 34795, 34797, 34799, 34803, 34805, 34807, 34811, 34813, 34815, 34943, 35007, 35039, 35055, 35063, 35067, 35069, 35071, 35135, 35167, 35183, 35191, 35195, 35197, 35199, 35231, 35247, 35255, 35259, 35261, 35263, 35279, 35287, 35291, 35293, 35295, 35303, 35307, 35309, 35311, 35315, 35317, 35319, 35323, 35325, 35327, 35391, 35423, 35439, 35447, 35451, 35453, 35455, 35487, 35503, 35511, 35515, 35517, 35519, 35535, 35543, 35547, 35549, 35551, 35559, 35563, 35565, 35567, 35571, 35573, 35575, 35579, 35581, 35583, 35615, 35631, 35639, 35643, 35645, 35647, 35663, 35671, 35675, 35677, 35679, 35687, 35691, 35693, 35695, 35699, 35701, 35703, 35707, 35709, 35711, 35727, 35735, 35739, 35741, 35743, 35751, 35755, 35757, 35759, 35763, 35765, 35767, 35771, 35773, 35775, 35783, 35787, 35789, 35791, 35795, 35797, 35799, 35803, 35805, 35807, 35815, 35819, 35821, 35823, 35827, 35829, 35831, 35835, 35837, 35839, 35903, 35935, 35951, 35959, 35963, 35965, 35967, 35999, 36015, 36023, 36027, 36029, 36031, 36047, 36055, 36059, 36061, 36063, 36071, 36075, 36077, 36079, 36083, 36085, 36087, 36091, 36093, 36095, 36127, 36143, 36151, 36155, 36157, 36159, 36175, 36183, 36187, 36189, 36191, 36199, 36203, 36205, 36207, 36211, 36213, 36215, 36219, 36221, 36223, 36239, 36247, 36251, 36253, 36255, 36263, 36267, 36269, 36271, 36275, 36277, 36279, 36283, 36285, 36287, 36295, 36299, 36301, 36303, 36307, 36309, 36311, 36315, 36317, 36319, 36327, 36331, 36333, 36335, 36339, 36341, 36343, 36347, 36349, 36351, 36383, 36399, 36407, 36411, 36413, 36415, 36431, 36439, 36443, 36445, 36447, 36455, 36459, 36461, 36463, 36467, 36469, 36471, 36475, 36477, 36479, 36495, 36503, 36507, 36509, 36511, 36519, 36523, 36525, 36527, 36531, 36533, 36535, 36539, 36541, 36543, 36551, 36555, 36557, 36559, 36563, 36565, 36567, 36571, 36573, 36575, 36583, 36587, 36589, 36591, 36595, 36597, 36599, 36603, 36605, 36607, 36623, 36631, 36635, 36637, 36639, 36647, 36651, 36653, 36655, 36659, 36661, 36663, 36667, 36669, 36671, 36679, 36683, 36685, 36687, 36691, 36693, 36695, 36699, 36701, 36703, 36711, 36715, 36717, 36719, 36723, 36725, 36727, 36731, 36733, 36735, 36751, 36759, 36763, 36765, 36767, 36775, 36779, 36781, 36783, 36787, 36789, 36791, 36795, 36797, 36799, 36807, 36811, 36813, 36815, 36819, 36821, 36823, 36827, 36829, 36831, 36839, 36843, 36845, 36847, 36851, 36853, 36855, 36859, 36861, 36863, 36991, 37055, 37087, 37103, 37111, 37115, 37117, 37119, 37183, 37215, 37231, 37239, 37243, 37245, 37247, 37279, 37295, 37303, 37307, 37309, 37311, 37327, 37335, 37339, 37341, 37343, 37351, 37355, 37357, 37359, 37363, 37365, 37367, 37371, 37373, 37375, 37439, 37471, 37487, 37495, 37499, 37501, 37503, 37535, 37551, 37559, 37563, 37565, 37567, 37583, 37591, 37595, 37597, 37599, 37607, 37611, 37613, 37615, 37619, 37621, 37623, 37627, 37629, 37631, 37663, 37679, 37687, 37691, 37693, 37695, 37711, 37719, 37723, 37725, 37727, 37735, 37739, 37741, 37743, 37747, 37749, 37751, 37755, 37757, 37759, 37775, 37783, 37787, 37789, 37791, 37799, 37803, 37805, 37807, 37811, 37813, 37815, 37819, 37821, 37823, 37831, 37835, 37837, 37839, 37843, 37845, 37847, 37851, 37853, 37855, 37863, 37867, 37869, 37871, 37875, 37877, 37879, 37883, 37885, 37887, 37951, 37983, 37999, 38007, 38011, 38013, 38015, 38047, 38063, 38071, 38075, 38077, 38079, 38095, 38103, 38107, 38109, 38111, 38119, 38123, 38125, 38127, 38131, 38133, 38135, …